\documentclass[12pt,letterpaper,reqno]{amsart}

\usepackage{amsfonts}
\usepackage[normalem]{ulem}
\usepackage{xcolor}
\usepackage{bm} 
\usepackage{graphicx} 
\usepackage{tikz}
\usepackage{amsmath,amssymb,amsthm}
\usepackage{latexsym}
\usepackage[utf8x]{inputenc}
\usepackage[T1]{fontenc}
\usepackage{fancyhdr}

\renewcommand{\S}{\mathbb{S}}
\newcommand{\R}{\mathbb{R}}

\usepackage[colorlinks=true,urlcolor=blue]{hyperref}

\theoremstyle{theorem}

\newtheorem{theorem}{Theorem}[section]

\title{The Spherical Gr\"unbaum Inequality}

\author[S. Myroshnychenko, D. Ryabogin, K. Tatarko, V. Yaskin] {Sergii Myroshnychenko, Dmitry Ryabogin, Kateryna Tatarko, and Vladyslav Yaskin}

\address{S.~Myroshnychenko, Department of Mathematics and Statistics,  University of the Fraser Valley, Abbotsford, BC V2S 7M8,  Canada}
	
\email{serhii.myroshnychenko@ufv.ca}

    \address{D.~Ryabogin, Department of Mathematics, Kent State University, Kent, OH 44242, USA}

    \email{ryabogin@math.kent.edu}

	\address{K.~Tatarko, Department of Pure Mathematics, University of Waterloo, Waterloo, ON N2L 3G1, Canada}
	
	\email{ktatarko@uwaterloo.ca}

	\address{V.~Yaskin, Department of Mathematical and Statistical Sciences, University of Alberta, Edmonton, AB T6G 2G1, Canada}
	
	\email{yaskin@ualberta.ca}

	\thanks{The first, third and fourth authors are supported in part by NSERC. The second  author is  supported in part by U.S.~National Science Foundation Grant DMS-2247771}

	\subjclass[2020]{Primary 52A20, 52A40}
	
	\keywords{Convex body, center of mass, Gr\"unbaum inequality, spherically convex set}

\begin{document}

\maketitle

\begin{abstract}
We prove an analogue of Gr\"unbaum's inequality on the sphere. 
Let $n \geq 3$ and let $K$ be a convex body on $\S^{n-1}\subset \mathbb R^n$ with centroid at $\theta\in \S^{n-1}$. Then for any $u\in \S^{n-1}$ that is orthogonal to $\theta$ we have $$\sigma(K\cap u^+) \ge \left(1-\frac{1}{n}\right)^{n-1} \sigma(K),$$
 where $\sigma$ denotes the spherical measure.
    The constant in this inequality is optimal.
\end{abstract}

\section{Introduction}

Let $K$ be a convex body in $\mathbb R^n$, that is, a convex compact set with non-empty interior. The center of mass (or centroid) of $K$ is the point in the interior of $K$ defined by 
$$c(K)=\frac{1}{\mathrm{vol}_n(K)} \int_K x \, dx,$$
where $\mathrm{vol}_n(K)$ is the volume (the $n$-dimensional Lebesgue measure) of $K$, and integration is taken with respect to Lebesgue measure. Translating $K$ if necessary, we may assume that its centroid is at the origin. 

 Denote by $\S^{n-1}$ the unit sphere in $\R^n$. For a vector $u\in \S^{n-1}$, define $u^\perp=\{x\in \mathbb R^n: \langle x,u\rangle =0\}$,  $u^+=\{x\in \mathbb R^n: \langle x,u\rangle \ge 0\}$, and $u^-=\{x\in \mathbb R^n: \langle x,u\rangle \le 0\}$.

The following result was obtained by Winternitz \cite[\textsection 23]{B} for $n=2$ and later extended to all dimensions by Gr\"unbaum \cite{G}. 
Let $K$ be a convex body in $\mathbb R^n$ with centroid at the origin. Then, for any $u\in \S^{n-1}$,
	\begin{equation}\label{Gru}
		\mathrm{vol}_n(K\cap u^+) \ge \left(1-\frac{1}{n+1}\right)^n \mathrm{vol}_n(K),
		\end{equation}
 and the inequality is sharp.

The goal of this paper is to prove an analogous result for convex bodies on the sphere $\S^{n-1}$.   Throughout the paper, we use $\sigma(A)$ to denote the spherical measure of a set $A\subset\mathbb S^{n-1}$.  
 We say that $A\subset\mathbb S^{n-1}$ is a convex body if it is contained in some open hemisphere and  $\mathcal C=\{ t \xi: \ 0\le t \le 1, \ \xi \in A \}$ is a convex body in $\R^n$. Then the centroid of $A$ is the radial projection onto $\S^{n-1}$ of the Euclidean centroid of $\mathcal C$ (see Section~\ref{Section::prelim} for details).
Our main result is the following.

\begin{theorem}\label{MainThm}
    Let $K$ be a convex body on $\S^{n-1}$, $n\ge 3$, with centroid at $
    \theta\in \S^{n-1}$. For any $u\in \S^{n-1}$ that is orthogonal to $\theta$ we have $$\sigma(K\cap u^+) \ge \left(1-\frac{1}{n}\right)^{n-1} \sigma(K).$$
    The constant in the latter inequality is optimal.
\end{theorem}
Some partial results in the case $n=3$ were obtained in \cite{BHLL}. Related problems have been studied  in spaces with non-negative Ricci curvature that contain an unbounded geodesic (see \cite{BOS}), and in the settings of $s$-concave measures \cite{MSZ} (see also \cite{FLLMT}). For more results on Gr\"unbaum-type inequalities, see \cite{AMMY}, \cite{BO}, \cite{CS}, \cite{FMY}, \cite{LY}, \cite {MNRY}, \cite{O} \cite{R}, \cite{ShY}, \cite{SY}, \cite{SZ},  and references therein.

\section{Preliminaries} \label{Section::prelim}
Let $\S^{n-1}_+=\mathbb S^{n-1}\cap \{x_n>0\}$ denote the open upper hemisphere of the unit sphere $\S^{n-1}\subset \mathbb R^n$, where $n\ge 3$. 
The  definition of spherical convexity given in the introduction is equivalent to the following. A set $A \subset \S^{n-1}_+$ is called convex if, for every pair of points $a, b \in A$, the geodesic segment that joins $a$ and $b$ lies in $A$. Since any two points in $\S^{n-1}_+$ can be connected by a unique geodesic segment (that is, an arc of a great circle), convexity is well defined in this setting.

We say that $A \subset \S^{n-1}_+$ is a spherically convex body if $A$ is closed, convex, and has non-empty interior relative to $\S^{n-1}$.
For a spherically convex body $A$, its  moment vector is $M(A)= \int_{A} \xi \,d\sigma(\xi)$, and its centroid is defined by 
$$c_s(A) = \frac{M(A)}{|M(A)|},$$
where $d\sigma$ denotes the  area measure on $\mathbb S^{n-1}$, and $|\cdot|$ is the Euclidean norm on $\mathbb R^n$. Alternatively, one can define $c_s(A)$ as the radial projection onto $\mathbb S^{n-1}$ of the Euclidean centroid of $\mathcal C=\{ t \xi: \ 0\le t \le 1, \ \xi \in A \}$, that is,
$$c_s(A) = \frac{c(\mathcal C)}{|c(\mathcal C)|}.$$

The gnomonic projection of $ \mathbb S^{n-1}_+$ onto the hyperplane $\{x_n=1\}$, which is the tangent hyperplane to the sphere $\mathbb S^{n-1}$ at the north pole $N=(0,\ldots,0,1)$, is defined as follows. For a point $q\in  \mathbb S^{n-1}_+$, its gnomonic projection is the point $\tilde q$ obtained as the point of intersection of the hyperplane $\{x_n=1\}$ and the line passing through the origin and the point $q$. For a convex body $A\subset  \mathbb S^{n-1}_+ $, its gnomonic projection $\tilde A$ is a convex body in $\{x_n=1\}$. Moreover, the following relations hold:
 \begin{equation}\label{measure-gnomonic}
     \sigma(A) = \int_{\tilde A} (1+|x|^2)^{-\frac{n}{2}} \, dx_1\ldots dx_{n-1},
     \end{equation}
 and the gnomonic projection of $c_s(A)$ is
\begin{equation}\label{proj-of-centroid} \tilde c_s(A) =\frac{\int_{\tilde{A}} x(1+|x|^2)^{-\frac{n+1}{2}} \, dx_1\ldots dx_{n-1}}{\int_{\tilde{A}} (1+|x|^2)^{-\frac{n+1}{2}} \, dx_1\ldots dx_{n-1} };
\end{equation}
 see \cite[\textsection 2]{BHPS}.

If $p$ and $q$ are two points in $\mathbb R^n$, the straight line segment connecting these two points will be denoted by $[p,q]$. If $p$ and $q$
are points in $\S^{n-1}_+$, we will write $[p,q]_{\S}$ for the geodesic arc connecting these points. The spherical distance between the points $p$ and $q$
 in $\S^{n-1}_+$ is the length of the geodesic arc $[p,q]_{\S}$, and will be denoted by $\mathrm{dist}_{\S}(p,q)$. We say that  the hyperplane $u^\perp$ intersects a  spherical arc $[p,q]_{\S}$ orthogonally, if  the tangent vector to $[p,q]_{\S}$  at the point of its intersection with $u^\perp$ is parallel to $u$.

We say that a set $\Omega\subset \mathbb R^n$ is star-shaped with respect to the origin if for every point $q\in \Omega$, any point in the interval between the origin and $q$ lies in $\Omega$.
For a star-shaped set $\Omega\subset \mathbb R^n$, its radial function is given by
 $$\rho_\Omega(\xi) = \sup\{t\ge 0: t\xi\in \Omega\}, \qquad \xi \in \S^{n-1}.$$
 The volume of $\Omega$ in polar coordinates equals
 $$\mathrm{vol}_n(\Omega)=\frac1n \int_{\S^{n-1}} \rho_{\Omega}^n(\xi) \, d\sigma(\xi).$$
 The Minkowski sum of two sets $\Omega_1, \Omega_2 \subset \R^n$ is defined as $\Omega_1 + \Omega_2 = \{q_1 + q_2 : q_1 \in \Omega_1, q_2 \in \Omega_2\}$.
The Hausdorff distance between two compact sets $\Omega_1, \Omega_2$ in $\mathbb{R}^n$ is 
$$
\delta(\Omega_1, \Omega_2) = \min \{ \lambda \geq 0: \, \Omega_1 \subset \Omega_2 + \lambda B^n_2, \, \Omega_2 \subset \Omega_1 + \lambda B^n_2 \},
$$
where $B_2^n$ is the unit Euclidean ball in $\mathbb{R}^n$  centered at the origin.

We will also need the following generalization of Gr\"unbaum’s inequality for sections of convex bodies:
\begin{theorem}[\cite{MSZ}]\label{Grunbaum_sections}
    Let $E$ be a $k$-dimensional subspace of $\R^n$ and let $\eta \in  \S^{n-1}\setminus E^\perp$. 
    Let $K$ be a convex body in $\R^n$ with  centroid  $c(K) \in E \cap \eta^\perp$. Then
    $$
    \frac{\mathrm{vol}_k(K \cap E \cap \eta^+)}{\mathrm{vol}_k(K \cap E)} \geq  \left(\frac{k}{n + 1}\right)^{k}.
    $$
    \end{theorem}

 In particular, when $E$ is 1-dimensional, this yields the well-known result due to Minkowski and Radon \cite[p.~58]{BF}:
    \begin{equation}\label{MR}
    \frac{\mathrm{vol}_1(K \cap E \cap \eta^+)}{\mathrm{vol}_1(K \cap E)} \geq   \frac{1}{n + 1} .
    \end{equation}

\section{Proof of Theorem \ref{MainThm}}

For the reader's convenience, we begin with an outline of the proof. Let $u \in \S^{n-1}$ and $K$ be a convex body on $\S^{n-1}$ with centroid $c_s(K) \in u^\perp$.

\begin{figure}[h!]
	\includegraphics[scale=0.2]{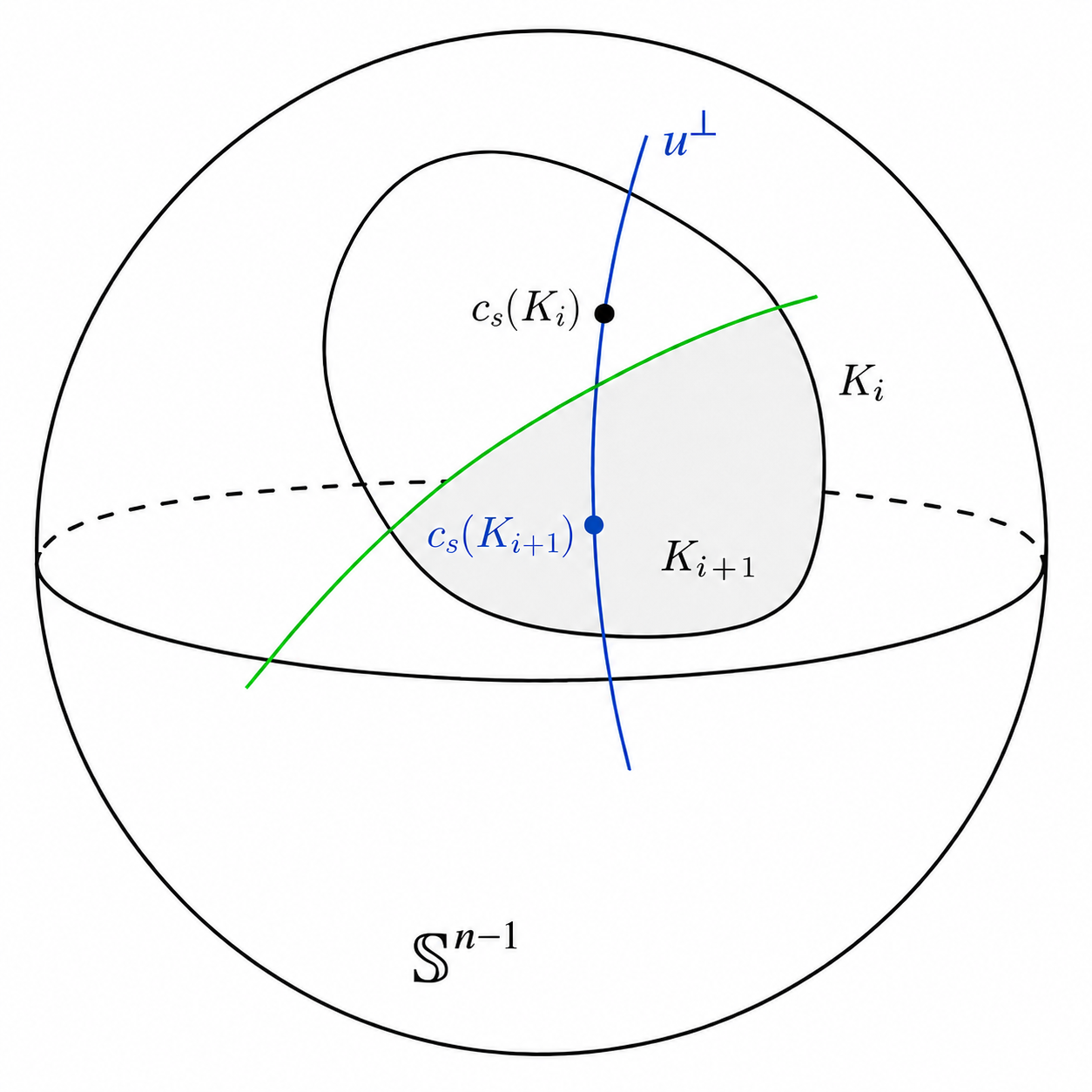}
	\caption{ Construction of $K_{i+1}$ from $K_i$ in Step 1.}
	\label{Fig:slic}
\end{figure}

{\bf Step 1}. First, we produce a sequence of spherically convex bodies $\{K_i\}_{i=0}^\infty$, such that
\begin{itemize}
    \item  $K_0= K$, 
    \item each $K_{i+1}$ is obtained from $K_i$ by intersecting $K_i$ with a certain half-space, so that $K_{i+1}\subset K_i$,

    \item $c_s(K_i)\in u^\perp$, for each $i\ge 0$,
    \item $$			\frac{\sigma(K_i \cap u^+)}{\sigma(K_i )}\ge \frac{\sigma(K_{i+1} \cap u^+)}{\sigma(K_{i+1} )}
			$$
            for each $i\ge 0$ (see Figure~\ref{Fig:slic}).

\end{itemize}
            
The sequence $\{K_i\}_{i=0}^\infty$ converges to 
		$K_\infty=  \bigcap_{i= 0}^\infty K_i$ in the Hausdorff metric, such that $K_\infty\cap u^\perp$ is a single point. Then we show that $K_\infty$ is either a point or a spherical segment.

        Our goal is to prove that
		\begin{equation}\label{eq::main}
		    \lim_{i\to \infty} \frac{\sigma(K_i \cap u^+)}{\sigma(K_i )} \ge \left(1-\frac{1}{n}\right)^{n-1}.
		\end{equation}

    In Steps 2 - 5, we focus on the case where $K_\infty$ is a spherical segment,  and in Step 6 we deal with the case when $K_\infty$ is a point.

{\bf Step 2.} We replace $u$ with a vector $v$ such that
    
$$ K_\infty\cap v^\perp=K_\infty\cap u^\perp=\lim_{i\to \infty}c_s(K_i)$$ and  $v^\perp$ intersects $K_\infty$ orthogonally.

        For this $v$ we obtain
        $$\lim_{i\to \infty} \frac{\sigma(K_i \cap u^+)}{\sigma(K_i )}=\lim_{i\to \infty} \frac{\sigma(K_i \cap v^+)}{\sigma(K_i )}.$$

{\bf Step 3.} We  pass to a sequence of bodies $\{\Omega_i\}$ in $\mathbb R^n$  defined as
$$\Omega_i=\{t \xi: t \ge 0, \xi\in K_i, t( \xi_1^2+\xi_n^2)\le 1\}.$$ 
For this sequence, we have 
$$\lim_{i\to \infty} \frac{\sigma(K_i \cap v^+)}{\sigma(K_i )}=\lim_{i\to \infty}\frac{\mathrm{vol}_n (\Omega_i\cap v^+)}{\mathrm{vol}_n (\Omega_i)}.$$
 Thus, to prove \eqref{eq::main}, it suffices to show that 
$$
\lim_{i\to \infty}\frac{\mathrm{vol}_n (\Omega_i\cap v^+)}{\mathrm{vol}_n (\Omega_i)} \ge \left(1-\frac{1}{n}\right)^{n-1}.
$$

{\bf Step 4.} 
The hyperplane $v^\perp$ does not necessarily contain the centroids of the bodies $\Omega_i$. 
We show that we can replace $v$ with $v_i$ such that $v_i^\perp$ is orthogonal to $K_\infty$ and $c(\Omega_i)\in v_i^\perp$, for each $i$. 
Moreover, 
$$\lim_{i\to \infty}\frac{\mathrm{vol}_n (\Omega_i\cap v^+)}{\mathrm{vol}_n (\Omega_i)} =\lim_{i\to \infty}\frac{\mathrm{vol}_n (\Omega_i\cap v_i^+)}{\mathrm{vol}_n (\Omega_i)}. $$

{\bf Step 5.} We transform $\Omega_i$ into a cone $\hat \Omega_i$ that is the convex hull of the origin and an $(n-1)$-dimensional convex set $A_i$  which lies in an $(n-1)$-dimensional affine plane not passing through the origin, such that
\begin{equation*}\frac{\mathrm{vol}_n (\Omega_i\cap v_i^+)}{\mathrm{vol}_n (\Omega_i)} = \frac{\mathrm{vol}_n (\hat\Omega_i\cap v_i^+)}{\mathrm{vol}_n (\hat\Omega_i)}.
        \end{equation*}
 The centroid of $\hat \Omega_i$ may have shifted to $v_i^+$. We find a vector $w_i$ such that $c(\hat \Omega_i)\in w_i^\perp$ and $w_i^\perp$ is orthogonal to $K_\infty$. Additionally,
$$
\frac{\mathrm{vol}_n (\hat\Omega_i\cap v_i^+)}{\mathrm{vol}_n (\hat\Omega_i)}
        \ge \frac{\mathrm{vol}_n (\hat\Omega_i\cap w_i^+)}{\mathrm{vol}_n (\hat\Omega_i)} =\frac{\mathrm{vol}_{n-1} (A_i\cap w_i^+)}{\mathrm{vol}_{n-1} (A_i)}.
        $$

For the convex body $A_i$ in $\mathbb R^{n-1}$ with centroid in $w_i^\perp$, we use \eqref{Gru} to conclude 
$$
\frac{\mathrm{vol}_{n-1} (A_i\cap w_i^+)}{\mathrm{vol}_{n-1} (A_i)}\ge \left(1-\frac{1}{n}\right)^{n-1}.
$$

Hence, the proof of \eqref{eq::main} is complete in the case where $K_\infty$ is a segment.

{\bf Step 6}. We now consider the case where $K_\infty$ is a single point. 
For sufficiently large $i$, the sets $K_i$ lie in small spherical neighbourhoods of $K_\infty$ where the spherical geometry is approximately Euclidean. This allows us to apply \eqref{Gru} in the limit.

{\bf Step 7.} We prove that the constant $\left(1-\frac{1}{n}\right)^{n-1}$ in 
Theorem \ref{MainThm} is optimal by taking a sequence of spherical simplices converging to a point.

\bigskip

We now provide more details for each step.

\bigskip

{\bf Step 1.}
Let $K$ be a spherically convex body on $ \mathbb{S}^{n-1}$, $n \geq 3$. Without loss of generality, we can assume that the vector $u$ from the statement of the theorem has the form $u=\cos \alpha e_{n-1}+\sin\alpha e_n$, for some angle $\alpha \in (-\pi\slash2, \pi\slash2)$. Here and below, $\{e_1, \ldots, e_n\}$ denotes the standard orthonormal basis in $\mathbb R^n$. Recall that $K\subset \mathbb S^{n-1}_+$ and  $c_s(K)\in u^\perp$.   
				
We  define recursively a sequence of spherically convex bodies $\{K_i\}_{i=0}^\infty$ as follows. Let $K_0=K$.  Suppose we have already constructed $K_i$ for some $i\ge 0$. Let $L_i=K_i\cap u^\perp$, and let $\tilde L_i$ be the gnomonic projection of $L_i$. 		Note that $\tilde L_i$ is contained in the $(n-2)$-dimensional affine subspace $Y= \{x_n=1\}\cap u^\perp$ that is orthogonal to $e_{n-1}$ and $e_n$. 
Consider the smallest $(n-2)$-dimensional box $B_i$ in $Y$, containing $\tilde L_i$, whose facets are parallel to the coordinate hyperplanes $e_1^\perp, \ldots, e_{n-2}^\perp$. Assume that the box is given by $$B_i= \{ x\in \mathbb R^n: a_1^{(i)}\le x_1\le b_1^{(i)} ,   \ldots, a_{n-2}^{(i)}\le x_{n-2}\le b_{n-2}^{(i)}, x_{n-1}= -\tan \alpha, x_n=1 \}.$$ 
	
Suppose that the longest edge of $B_i$ is parallel to the axis $e_{k_i}$ for some $1\le k_i\le n-2$ (if such an edge is not unique, choose any one of them). 			
Consider the $(n-3)$-dimensional affine subspace  $$V_i=\{ x\in \mathbb R^n:   x_{k_i}=(a_{k_i}^{(i)}+ b_{k_i}^{(i)})/2, x_{n-1}= -\tan\alpha, x_n=1 \}$$ contained in $Y$, which divides the box $B_i$ in half. Note that $V_i$ splits $\tilde L_i$ in two convex subsets of positive $(n-2)$-dimensional measure.
Let $H_i$ be  the linear span of $V_i$ in $\mathbb R^n$. We claim that there is a vector $\xi_i\in \S^{n-1}\cap H_i^\perp$ such that $\xi_i^\perp$ divides $K_i$  into two subsets $K_{i}^-=K_i\cap \xi_i^-$ and $K_i^+=K_i\cap \xi_i^+$, with $c_s(K_i^-)$ and $c_s(K_i^+)$ both lying in $u^\perp$.

      To prove the claim, it is enough to guarantee that the moment vectors   $M(K_i^{\pm}) =  \int_{K_i^{\pm}} \xi \,d\sigma(\xi)$ lie in $u^\perp$.
		
        Let $\langle u \rangle=\{x\in \mathbb R^n:x=t u, t\in \mathbb R\}$.  Let $f: \S^{n-1}\cap H_i^\perp \to \langle u \rangle$ be defined by
        $$f(\xi)= P_{\langle u \rangle} M(K_i \cap \xi^+), 
        $$
        where $P_{\langle u \rangle}$ is the orthogonal projection onto the line $\langle u \rangle$.   We note that we can identify $\S^{n-1}\cap H_i^\perp$ with $\S^{1}$ and $\langle u \rangle$ with $\R.$
        
        Observe that, for all $\xi\in \S^{n-1} \cap H_i^\perp$, 
		\begin{align*}
			f(-\xi) &= P_{\langle u \rangle} M(K_i \cap (-\xi)^+) = P_{\langle u \rangle} M(K_i \cap {\xi}^-) = P_{\langle u \rangle} M(K_i \backslash (K_i \cap {\xi}^+)) \\
			&= P_{\langle u \rangle} M  (K_i) - P_{\langle u \rangle} M(K_i \cap {\xi}^+) = - P_{\langle u \rangle} M(K_i \cap {\xi}^+) = -f(\xi),
		\end{align*}
		where $P_{\langle u \rangle} M  (K_i) = 0$ because $c_s(K_i) \in u^\perp$. Since $f$ is continuous and odd, we conclude that there is $\xi_i \in  \S^{n-1}\cap H_i^\perp $ such that $f(\xi_i) = 0$. Thus, $M(K_i \cap {\xi_i}^+)\in u^\perp$, and therefore, $M(K_i \cap {\xi_i}^-)\in u^\perp$ since $c_s(K_i)$ is a linear combination of $c_s(K_i^-)$ and $c_s(K_i^+)$, as claimed.

	  Now we define the set $K_{i+1}$. Note that $\sigma(K_{i}^-)>0$ and $\sigma(K_{i}^+)>0$  since $K_{i}^-$ and $K_{i}^+$ are sets obtained by cutting $K_i$ by a hyperplane passing through the relative interior of $K_i$.  Without loss of generality, we can assume that 
		$$
		\frac{\sigma(K_{i}^-\cap u^+)}{\sigma(K_{i}^-)}\le	\frac{\sigma(K_{i}^+\cap u^+)}{\sigma(K_{i}^+)}.$$
		A simple computation shows that the latter  implies
		$$
		\frac{\sigma(K_i^-\cap u^+)}{\sigma(K_i^-)}\le \frac{\sigma(K_i^-\cap u^+)+\sigma(K_i^+\cap u^+)}{\sigma(K_i^-)+\sigma(K_i^+)}= \frac{\sigma(K_i \cap u^+)}{\sigma(K_i)}.$$
		Denoting $K_{i+1}=K_i^-$, we obtain a spherically convex body $K_{i+1}$ with $c_s(K_{i+1})\in u^\perp$ such that 
		\begin{equation}\label{decrease}
			\frac{\sigma(K_{i+1} \cap u^+)}{\sigma(K_{i+1} )}\le \frac{\sigma(K_i \cap u^+)}{\sigma(K_i )}.
			\end{equation}

		Thus, we have constructed a sequence of spherically convex bodies $\{K_i\}_{i=0}^\infty$ with  $c_s(K_i)\in u^\perp$ for all $i$, such that 
		 $K_{i+1}\subset K_{i}$, and  inequality \eqref{decrease} is satisfied (see Figure~\ref{Fig:slic}). Moreover, $K_i\cap u^\perp$ converges to a single point as $i\to \infty$, since $\tilde L_i \subset B_i$ for all $i$, and the longest edge of $B_i$ is divided in half at each step.

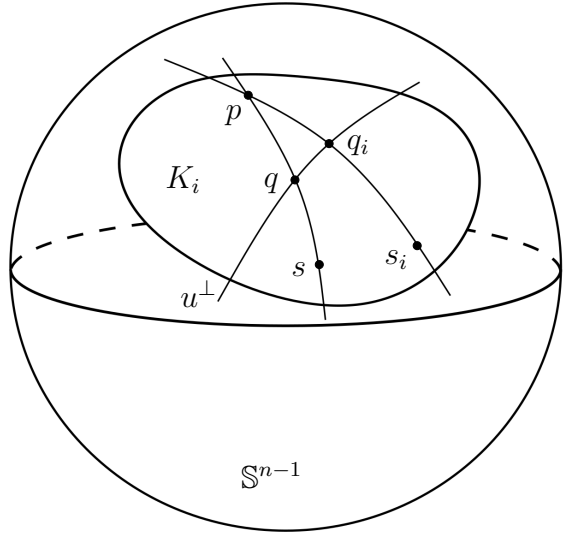
\begin{figure}[h!]
 \centering
\begin{tikzpicture}
\node[anchor=south west, inner sep=0] (img) at (0,0)
        {\includegraphics[scale=0.23]{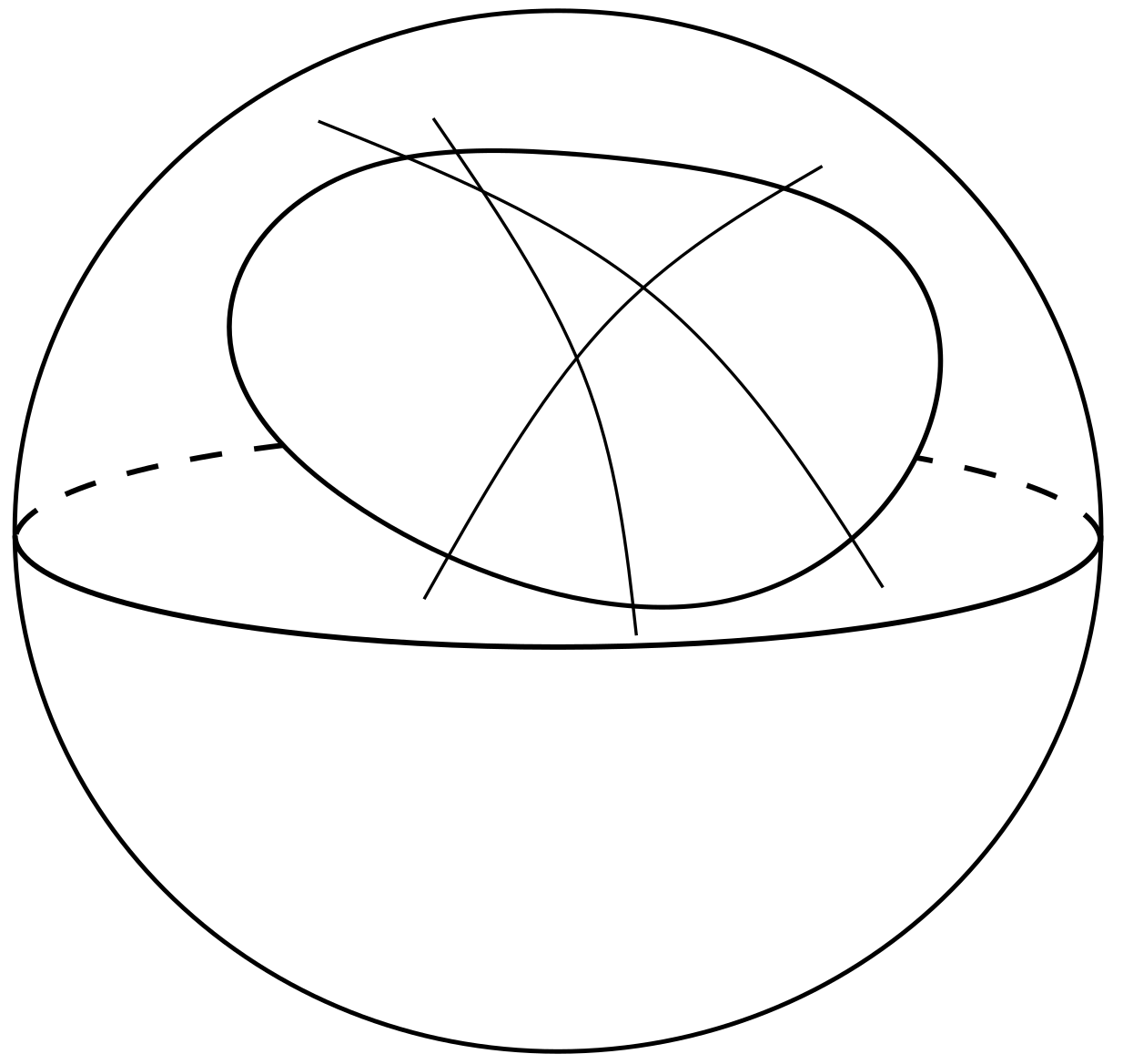}};
\draw [fill=black] (3.858,4.72) circle (1.5pt);
\draw [fill=black] (5.475,3.85) circle (1.5pt);
\draw [fill=black] (3.24,5.84) circle (1.5pt);
\draw [fill=black] (4.18,3.6) circle (1.5pt);
\draw [fill=black] (4.31,5.2) circle (1.5pt);
\node[above right] at (3.3,4.4) {$q$};
\node[above right] at (4.4,4.9) {$q_i$};
\node[above left] at (5.5,3.4) {$s_i$};
\node[above left] at (4.15,3.3) {$s$};
\node[above right] at (2.2,2.9) {$u^\perp$};
\node[above right] at (2.8,5.3) {$p$};
\node[above right] at (3,0.5) {$\mathbb{S}^{n-1}$};
\node[above right] at (2,4.4) {$K_i$};
\end{tikzpicture}
\caption{Constructing arc $[s,p]_{\S}$ in $K_\infty$.}
	\label{Fig:conv}
\end{figure}

		Let $$K_\infty=  \bigcap_{i= 0}^\infty K_i.$$
		Since  $K_i\cap u^\perp$ converges to a point $q$ as $i\to \infty$, we claim that $K_\infty$ can only be a spherical segment or a point.   To prove the claim, assume that $K_\infty$ is not a point. Then there exists a point $p\in K_\infty$ that is different from $q=K_\infty\cap u^\perp$. For each $i$, let $q_i\in K_i\cap u^\perp$ be the centroid of $K_i$, and let $E_i$ be the linear span of $p$ and $q_i$.
        Since the spherical segment $[p,q_i]_{\S}$ lies in $K_i$, we can use Theorem \ref{Grunbaum_sections} for the body $\mathcal C_i=\{ t \xi: \ 0\le t \le 1, \ \xi \in K_i \}$ with $E=E_i$  and $\eta=   u$ or $\eta=-  u$, chosen so that $\eta^+$ does not contain $p$.   Then we conclude that, for each $i$,
		\begin{align*}
		\mathrm{length}(K_i \cap E_i \cap \eta^+)&= 2\, {\mathrm{vol}_2(\mathcal C_i \cap E_i \cap \eta^+)}
        \geq  2 \left(\frac{2}{n + 1}\right)^{2} {\mathrm{vol}_2(\mathcal C_i \cap E_i)}\\
        &= \left(\frac{2}{n + 1}\right)^{2} \mathrm{length}(K_i \cap E_i) \ge  \left(\frac{2}{n + 1}\right)^{2} \mbox{dist}_{\mathbb S} (p,q_i).
		\end{align*}
		Hence, there is a point $s_i$ that belongs to  $K_i$, where $p$ and $s_i$ lie in different open half-spaces defined by $u^\perp$, and such that $q_i\in [s_i,p]_{\S}$. Moreover, $$\mbox{dist}_{\S} (s_i,q_i)\ge  \left(\frac{2}{n + 1}\right)^{2} \mbox{dist}_{\S}(p,q_i).$$
		Passing to a subsequence if needed, we can assume that $s_i$ converges to a point $s\in K_\infty$, as $i\to \infty$. Also note that $q_i$ converges to $q$. Thus, the geodesic segment $[s,p]_{\S}$ lies in $K_\infty$, with $q$ strictly inside this segment (see Figure \ref{Fig:conv}).

		 Extending this  segment if needed, we can assume that $[s,p]_{\S}=K_\infty \cap g$, where $g$ is the great circle on $\mathbb S^{n-1}$ containing $s$ and $p$. To reach a contradiction, let us assume that $K_\infty\ne [s,p]_{\S}$. Then there is a point $z\in K_\infty$ that is not on $g$. Consider the spherical triangle with vertices $s$, $p$, and $z$. The triangle intersects $u^\perp$ in a segment, which is impossible, since $K_\infty\cap u^\perp$ is a single point.  Thus, $K_\infty$ is a spherical segment.

\begin{figure}[h!]
	\includegraphics[scale=0.5]{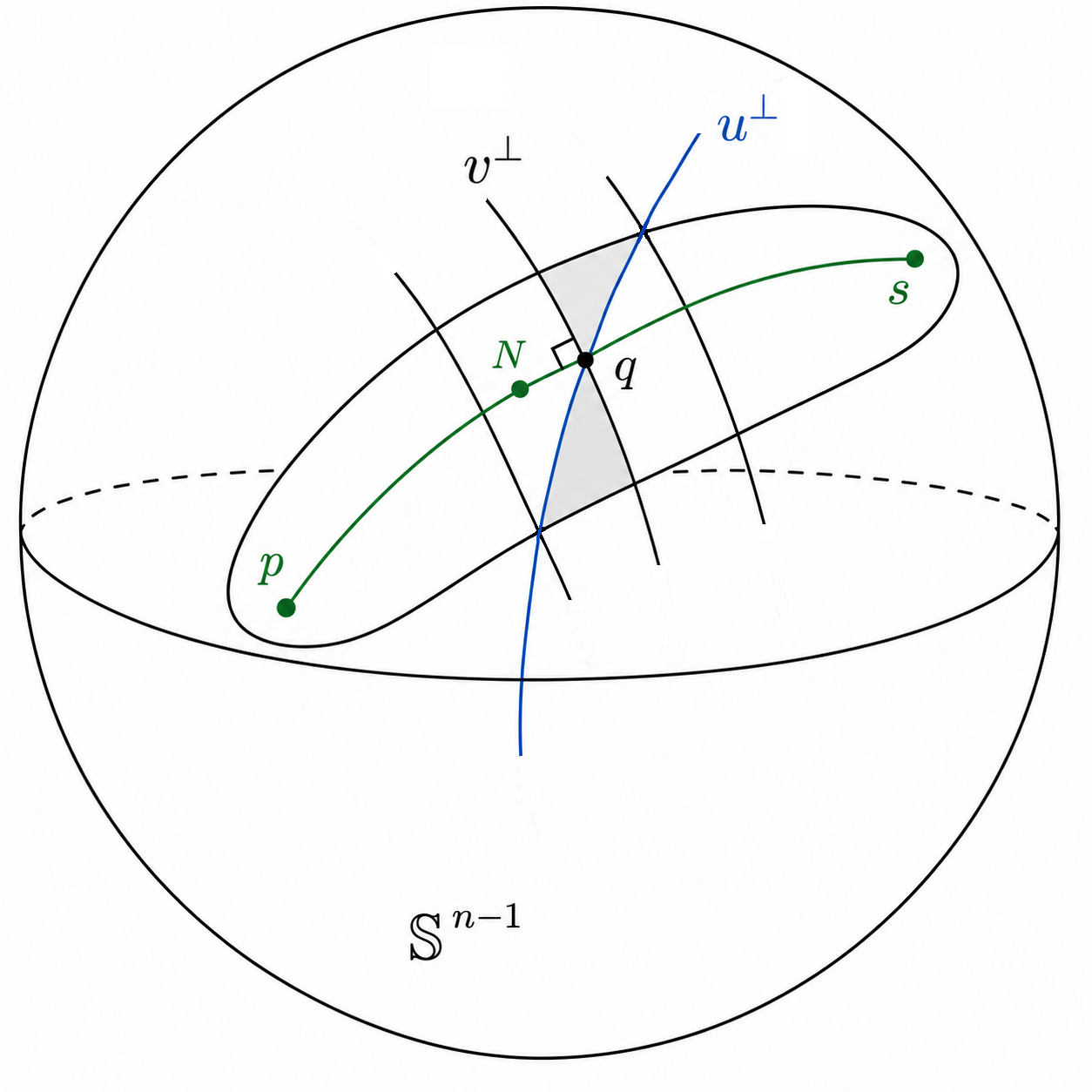}
	\caption{Replacing $u^\perp$ by $v^\perp$ in Step 2.}
	\label{Fig:turn}
\end{figure}

	{\bf Step 2.}	
		We first deal with the case where $K_\infty$ is a spherical segment. Let $s$ and $p$ be the end points of $K_\infty$, and let  $q= K_\infty \cap u^\perp$. Consider a unit vector $v$ such that $K_\infty \cap v^\perp=q$, the hyperplane $v^\perp$ intersects $K_\infty$ orthogonally, and $\langle p,v\rangle $ has the same sign as  $\langle p,u\rangle $. Generally, $u$ and $v$ are different vectors (see Figure~\ref{Fig:turn}). However, we show that for our purposes we can use $v$ instead of $u$, since 
        \begin{equation}\label{u-to-v}
            \lim_{i\to \infty} \frac{\sigma(K_i \cap u^+)}{\sigma(K_i )}=\lim_{i\to \infty} \frac{\sigma(K_i \cap v^+)}{\sigma(K_i )}.
        \end{equation}
    To prove the latter, it suffices to verify that
		$$\lim_{i\to \infty} \frac{\sigma(K_i \cap u^-\cap v^+)+\sigma(K_i \cap u^+\cap v^-)}{\sigma(K_i )} =0,$$
 where sets $K_i \cap u^-\cap v^+$ and $K_i \cap u^+\cap v^-$ are depicted as shaded regions in Figure~\ref{Fig:turn}.
	 
 Without loss of generality, we can assume that the midpoint of the segment $[s,p]_{\S}$ is located at the north pole $N=(0,...,0,1)$. 		 When $i$ is large enough, $K_i$ is contained in the upper hemisphere $\S^{n-1}_+$.  Consider the gnomonic projection centered at $N$. Let $\tilde K_i$ be the image of $K_i$ and let $[\tilde s, \tilde p]$ be the image of $[s,p]_{\S}$ under the gnomonic projection.  Applying a rotation if necessary, we can assume that $$[\tilde s, \tilde p]=\{ x\in \mathbb R^n: -a\le x_1\le a, x_2=0, ..., x_{n-1}=0, x_n=1\},$$
 for some $a>0$.
 
 Let $\tilde q= (b, 0,...,0,1)$ with $-a < b < a$.  Then $$v^\perp\cap \{x_n=1\} = \{x\in \mathbb R^n: x_1=b,x_2=0, ..., x_{n-1}=0, x_n=1\}.$$ 
Since $\{K_i\}$ converges to $[s,p]_{\S}$ in the Hausdorff metric as $i\to \infty$, we can choose a decreasing sequence of positive numbers $\{\varepsilon_i\}$ that converges to zero such that $\tilde K_i$ is contained in the set $[\tilde s, \tilde p]+\varepsilon_i B_2^{n-1}$ for every $i$; here $ B_2^{n-1}$ is the unit Euclidean ball in $\{x_n=1\}$, and the sum is the Minkowski sum.

		 Observe that the image of the sets $K_i \cap u^-\cap v^+$ and $K_i \cap u^+\cap v^-$ under the gnomonic projection is contained in the slab 
		 $$S_i=\{ x\in \mathbb R^n: b- \varepsilon_i\tan \gamma \le x_1\le b+ \varepsilon_i\tan \gamma  ,   x_n=1\},$$ where $\gamma$ is the angle between $u^\perp \cap \{x_n=1\}$ and $v^\perp \cap \{x_n=1\}$ as hyperplanes in $\{x_n=1\}$. 
		
	 Using \eqref{measure-gnomonic}, we obtain 
	\begin{align*}\sigma(K_i \cap u^-\cap v^+) &= \int_{\tilde   K_i \cap u^-\cap v^+} (1+|x|^2)^{-n/2} dx_1\ldots dx_{n-1}\\
		&\le \int_{\tilde   K_i \cap u^-\cap v^+}  dx_1\ldots dx_{n-1}\\
		&\le  \int_{\tilde   K_i \cap S_i}  dx_1\ldots dx_{n-1} \\
		& = \int_{b- \varepsilon_i\tan \gamma}^{b+ \varepsilon_i\tan \gamma}\mathrm{vol}_{n-2}( \tilde  K_i \cap \{x_1=t\})  dt\\
	&	\le   2 \varepsilon_i\tan \gamma \, \max_{t}\mathrm{vol}_{n-2}( \tilde  K_i \cap \{x_1=t\})  .
		\end{align*}
		Now we show that $\max_{t}\mathrm{vol}_{n-2}( \tilde  K_i \cap \{x_1=t\})$ is bounded above by a constant multiple of $\sigma(K_i)$. If this maximum is attained for some $t_0\le 0$, let $G$ be the cone in $\{x_n=1\}$ with base $ \tilde  K_i \cap \{x_1=t_0\}$ and apex $(a,0,...,0,1)$. If  $t_0\ge 0$, then let $G$ be the cone  with base $ \tilde  K_i \cap \{x_1=t_0\}$ and apex $(- a,0,...,0,1)$. In both cases $G$ is contained in $\tilde K_i$ and its height is at least $a$. Thus,
		$$\frac1{n-1} a \,\mathrm{vol}_{n-2}( \tilde  K_i \cap \{x_1=t_0\})\le \mathrm{vol}_{n-1}(G)\le \mathrm{vol}_{n-1}(\tilde K_i).$$
		Hence,
		$$\sigma(K_i \cap u^-\cap v^+)\le   2 \varepsilon_i\tan \gamma \frac{n-1}{a} \mathrm{vol}_{n-1}(\tilde K_i).$$
       Since $|x_1| \leq a + \varepsilon_i$ and $x_2^2 + \dots + x_{n-1}^2 \leq \varepsilon_i^2$ for every $x = (x_1, \dots, x_{n-1}) \in \tilde K_i$, we obtain
		$$ (1+(a+\varepsilon_i)^2+\varepsilon_i^2)^{-n/2} \mathrm{vol}_{n-1}(\tilde K_i)\le  \int_{\tilde   K_i } (1+|x|^2)^{-n/2} dx_1\ldots dx_{n-1}=  \sigma(K_i).$$
		Therefore,
		$$\sigma(K_i \cap u^-\cap v^+)\le D \varepsilon_i \sigma(K_i)$$
		for some absolute constant $D>0$.
		
		Similarly, we obtain
			$$\sigma(K_i \cap u^+\cap v^-)\le D \varepsilon_i \sigma(K_i),$$ as required.
		Thus,  \eqref{u-to-v} holds.
        
        We remark that the distance from any point of $K_i \cap u^\perp$ to $q$ is of order $O(\varepsilon_i)$.
        Since $q$ belongs to $v^\perp$, and $c_s(K_i)$ lies in $K_i\cap u^\perp$, we conclude that the distance from $c_s(K_i)$ to $v^\perp$ is of order $O(\varepsilon_i)$.

		{\bf Step 3.} Let $$\Omega_i=\{t \xi: t \ge 0, \xi\in K_i, t( \xi_1^2+\xi_n^2)\le 1\}$$
		be the intersection of the solid cylinder $\{x\in \mathbb R^n: x_1^2+x_n^2\le 1\}$ and the infinite convex  cone  $U_i = \{t \xi: t \ge 0, \xi\in K_i \}$ in $\R^n$.
        We now show that 
        $$
        \lim_{i\to \infty} \frac{\sigma(K_i \cap v^+)}{\sigma(K_i )}=\lim_{i\to \infty}\frac{\mathrm{vol}_n (\Omega_i\cap v^+)}{\mathrm{vol}_n (\Omega_i)}.
        $$

        Let $\rho_{\Omega_i}$ be the radial function of $\Omega_i$. We note that $\rho_{\Omega_i}(\xi)= \frac{1}{\xi_1^2 + \xi_n^2}$ for all $\xi \in K_i$. Since the Hausdorff distance between $K_i$ and  the  arc $[s,p]_{\S}$ is $O(\varepsilon_i)$, the quantity $\xi_2^2 + \dots + \xi_{n-1}^2$ is also $O(\varepsilon_i)$ and thus $\xi_1^2 + \xi_n^2 = 1 - O(\varepsilon_i)$ for all $\xi \in K_i$. Hence, $\rho_{\Omega_i}(\xi)= 1 + O(\varepsilon_i)$ for $\xi \in K_i$.

		Therefore, passing to polar coordinates, we obtain
		$$\mathrm{vol}_n (\Omega_i) = \frac{1}{n} \int_{K_i} \rho_{\Omega_i}^n(\xi) \, d\sigma(\xi)=\frac{1}{n} \int_{K_i}(1+O(\varepsilon_i))\, d\sigma(\xi)=\frac{1}{n}(1+O(\varepsilon_i)) \sigma(K_i).$$
		Similarly,
		$$\mathrm{vol}_n (\Omega_i\cap v^+) = \frac{1}{n} \int_{K_i\cap v^+} \rho_{\Omega_i}^n(\xi) \, d\sigma(\xi) =\frac{1}{n}(1+O(\varepsilon_i)) \sigma(K_i\cap v^+).$$
		
		Therefore,
		
		\begin{equation}\label{K-to-Omega}\frac{\mathrm{vol}_n (\Omega_i\cap v^+)}{\mathrm{vol}_n (\Omega_i)}=\frac{\sigma(K_i \cap v^+)}{\sigma(K_i )} +O(\varepsilon_i).
        \end{equation}

		{\bf Step 4.}
		Note that $v^\perp$ does not necessarily pass through the centroid of $\Omega_i$, however a small rotation of $v^\perp$ does, as we  show now. Let $v=\cos \beta e_1+\sin \beta e_n$  for some angle $\beta \in (-\frac{\pi}2, \frac{\pi}2)$. Consider the vectors of the form $v_i=\cos \beta_i e_1+\sin \beta_i e_n$, for some angles $\beta_i$, so that $c(\Omega_i)\in v_i^\perp$, and the sign of $\langle p,v_i\rangle$ is the same as the sign of $\langle p,v\rangle$. Now we argue that $\beta - \beta_i = O(\varepsilon_i)$. Let $\mathcal C_i=\{ \lambda \xi: \ 0\le \lambda \le 1, \ \xi \in K_i \}$. The centroid of $\mathcal C_i$ lies on the same ray from the origin as the centroid of $K_i$, and    $c(\mathcal C_i)$ is closer to the origin than $c_s(K_i)$. 
        By the remark at the end of Step 2, the distance from $c_s(K_i)$ to $v^\perp$ is of order $O(\varepsilon_i)$, so the distance from  $c(\mathcal C_i)$ to $v^\perp$ is also of order $O(\varepsilon_i)$. Now observe that 
\begin{align*}\left|c(\Omega_i)-c(\mathcal C_i)\right|& = \left|\frac{1}{\mathrm{vol}_n (\Omega_i)} \int_{\Omega_i} x\, dx-\frac{1}{\mathrm{vol}_n (\mathcal C_i)} \int_{\mathcal C_i} x \,dx\right| \\
            & \le \frac{1}{\mathrm{vol}_n (\Omega_i)} \left|\int_{\Omega_i} x\, dx- \int_{\mathcal C_i} x \,dx \right| +\left|\frac{1}{\mathrm{vol}_n (  \Omega_i)}  -\frac{1}{\mathrm{vol}_n (\mathcal C_i)} \right|\left|\int_{\mathcal C_i} x \,dx\right|\\
             & \le \frac{1}{(n+1)\mathrm{vol}_n (\Omega_i)} \left|\int_{K_i} \xi \rho_{\Omega_i}^{n+1}(\xi) \,d\sigma(\xi)- \int_{K_i} \xi   \,d\sigma(\xi) \right|\\
             &\hspace{6cm}+\left|\frac{1}{\mathrm{vol}_n (  \Omega_i)}  -\frac{1}{\mathrm{vol}_n (\mathcal C_i)} \right|\mathrm{vol}_n (\mathcal C_i)\\
              & \le \frac{1}{(1+\frac 1 n)(1+O(\varepsilon_i)) \sigma(K_i)}  \int_{K_i}  \left( \rho_{\Omega_i}^{n+1}(\xi) -1\right)\,d\sigma(\xi)+\left|(1+O(\varepsilon_i))-1\right| \\
            &= O(\varepsilon_i).
			\end{align*}  
        Here, in the second inequality, we use the fact that $\rho_{\mathcal C_i}(\xi)=1$ for all $\xi \in K_i$.

		Thus, the distance between $c(\Omega_i)$ and $c(\mathcal C_i)$ is of order $O(\varepsilon_i)$, and therefore the distance from $c(\Omega_i)$ to $v^\perp$ is also of order $O(\varepsilon_i)$.

         Since $c(\Omega_i)$ lies in $v_i^\perp$, it can be represented as follows:
        $$c(\Omega_i)=\chi_1\zeta_1+\chi_2\zeta_2,$$
        where $\chi_1>0$ and $\chi_2$ are real numbers,  $\zeta_1= -\sin \beta_i e_1+\cos \beta_i e_n$, and $\zeta_2$ is a unit vector orthogonal to the linear span of $e_1$ and $e_n$.
        Since the Hausdorff distance between $K_i$ and  the  arc $[s,p]_{\S}$ is $O(\varepsilon_i)$, and the arc lies in the $x_1x_n$-plane, the coefficient $\chi_2$ is of order 
        $O(\varepsilon_i)$. By \eqref{MR}, the distance from $c(\Omega_i)$ to the origin is at least $1/(n+1)$, and thus 
        $$\sqrt{\chi_1^2+\chi_2^2}\ge 1/(n+1), $$
        implying
        $$\chi_1\ge 1/(n+1)+ O(\varepsilon_i)\ge 1/(2n)$$
        for large enough $i$.

        Projecting $c(\Omega_i)$ onto $v$, we obtain
        $$O(\varepsilon_i)=\langle c(\Omega_i), v\rangle = \chi_1\langle \zeta_1, v\rangle = \chi_1\sin (\beta-\beta_i). $$
        Since $\chi_1$ is bounded from below by an absolute constant and $|\beta - \beta_i|$ is bounded away from $\pi$, we conclude that $\beta-\beta_i = O(\varepsilon_i).$

        We now proceed as in Step 2 to obtain 			\begin{equation}\label{v-to-vi}\frac{\mathrm{vol}_n (\Omega_i\cap v^+)}{\mathrm{vol}_n (\Omega_i)}= \frac{\mathrm{vol}_n (\Omega_i\cap v_i^+)}{\mathrm{vol}_n (\Omega_i)}+ O(\varepsilon_i).
        \end{equation}
                Indeed, since 
                $$\mathrm{vol}_n (\Omega_i)= \frac1{n} (1+O(\varepsilon_i)) \sigma(K_i),$$
                $$\mathrm{vol}_n (\Omega_i\cap v^+)= \frac1{n} (1+O(\varepsilon_i)) \sigma(K_i\cap v^+),$$
                $$\mathrm{vol}_n (\Omega_i\cap v_i^+)= \frac1{n} (1+O(\varepsilon_i)) \sigma(K_i\cap v_i^+),$$
                it is enough to show that 
                \begin{equation}\label{Bound}\left|\sigma(K_i\cap v^+)-\sigma(K_i\cap v_i^+)\right|=O(\varepsilon_i)\sigma(K_i).
                \end{equation}
                As in Step 2, let
                $$v^\perp\cap \{x_n=1\} = \{x\in \mathbb R^n: x_1=b,x_2=0, ..., x_{n-1}=0, x_n=1\},$$ 
                and let
                $$v_i^\perp\cap \{x_n=1\} = \{x\in \mathbb R^n: x_1=b_i,x_2=0, ..., x_{n-1}=0, x_n=1\},$$ 
                for some numbers $b_i\in (-a,a)$. Since $\beta-\beta_i$ is of order $O(\varepsilon_i)$, we have $|b-b_i|\le d \varepsilon_i$, for some absolute constant $d>0$. 
                Then, 
                \begin{align*}\left|\sigma(K_i\cap v^+)-\sigma(K_i\cap v_i^+)\right| & \le \int_{b- d\varepsilon_i }^{b+ d \varepsilon_i }\int_{\tilde  K_i \cap \{x_1=t\}}(1+|x|^2)^{-\frac{n}2} dx_2\ldots dx_{n-1} dt  \\              
                &\le \int_{b- d\varepsilon_i }^{b+ d \varepsilon_i } \mathrm{vol}_{n-2}( \tilde  K_i \cap \{x_1=t\})\,  dt\\
	&	\le   2d \varepsilon_i  \max_{t}\mathrm{vol}_{n-2}( \tilde  K_i \cap \{x_1=t\})  \\
    & \le  \frac{ 2d(n-1)}{a}  \varepsilon_i \mathrm{vol}_{n-1}(\tilde K_i)\\
    & \le  \frac{ 2d(n-1)}{a}  \varepsilon_i  (1+(a+\varepsilon_i)^2+\varepsilon_i^2)^{\frac{n}{2}} \sigma(K_i).
    \end{align*}
This proves \eqref{Bound}, and  \eqref{v-to-vi} follows.

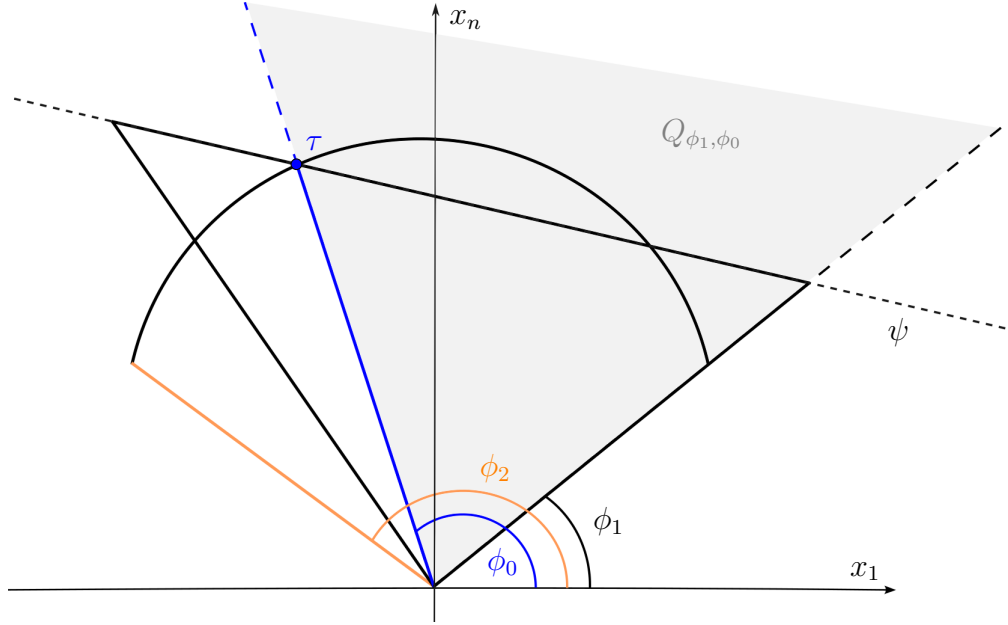
\begin{figure}[h!]
 \centering
\begin{tikzpicture}
\node[anchor=south west, inner sep=0] (img) at (0,0)
        {\includegraphics[scale=0.3]{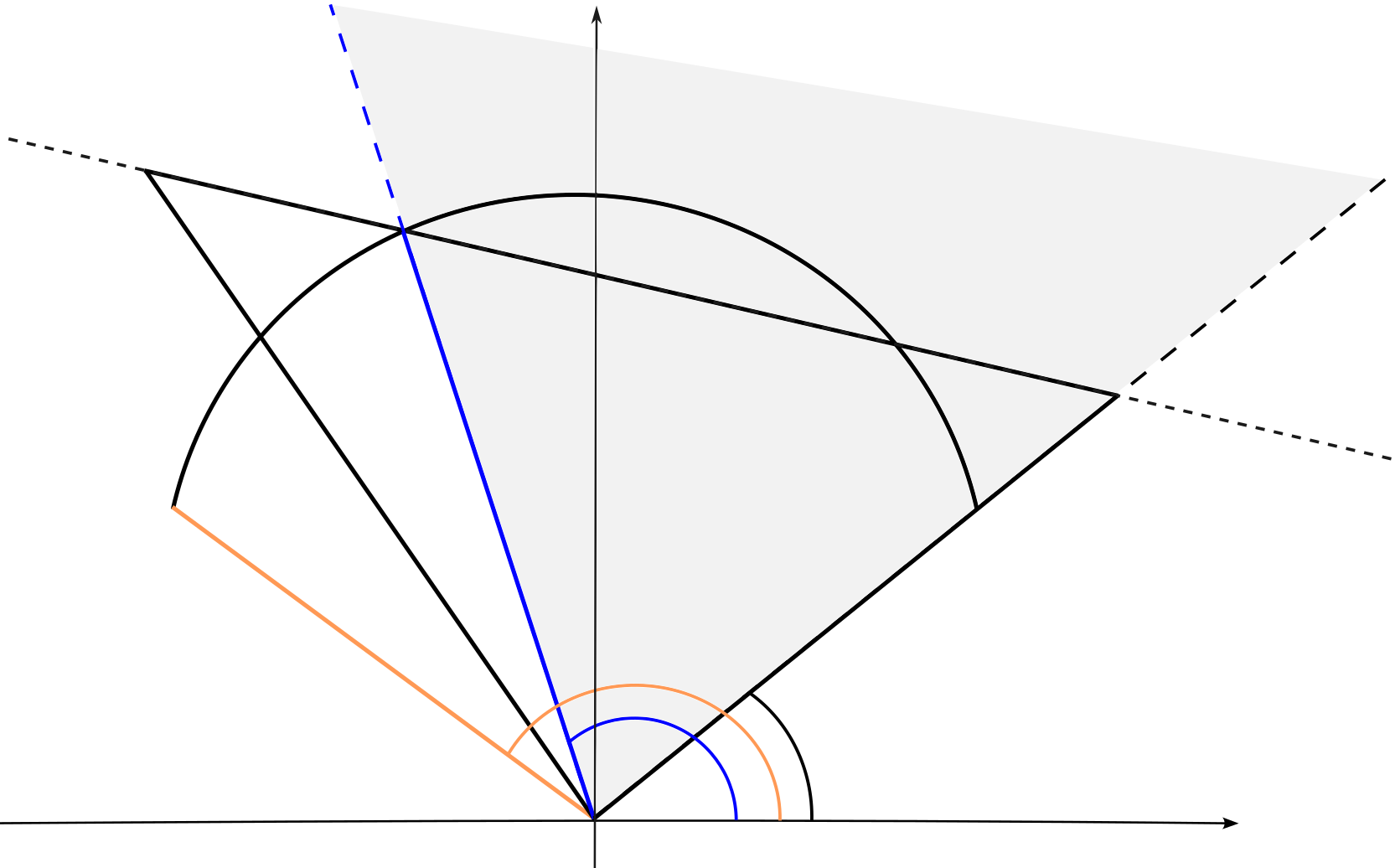}};
\draw [fill=blue] (3.82,6.05) circle (2pt);
\node[above right] at (3.8,6.1) {\textcolor{blue}{$\tau$}};
\node[above right] at (8.5,6.1) {\textcolor{gray}{$Q_{\phi_1,\phi_0}$}};
\node[above right] at (11.5,3.5) {$\psi$};
\node[above right] at (11,0.4) {$x_1$};
\node[above right] at (5.7,7.7) {$x_n$};
\node[above right] at (7.6,1) {$\phi_1$};
\node[above right] at (6.2,0.45) {\textcolor{blue}{$\phi_0$}};
\node[above right] at (6.1,1.7) {\textcolor{orange}{$\phi_2$}};
\end{tikzpicture}
\caption{Projection of $\Omega_i$ onto the $x_1x_n$-plane. }
	\label{Fig:icecream-corners}
\end{figure}

		{\bf Step 5.} Recall that $v_i = \cos \beta_i e_1+\sin \beta_i e_n$ for $\beta_i \in \left(-\frac \pi 2, \frac \pi 2 \right).$ Let $W_i=\{x\in \mathbb R^n: x_1^2+x_n^2= 1\} \cap v_i^\perp$ and let $\Psi_i$ be the  one-parametric family of all affine hyperplanes containing $W_i$. 		
		Suppose that the projection of $\Omega_i$ onto the two-dimensional coordinate plane $x_1x_n$ is the sector that in polar coordinates ($x_1=r\cos\phi$, $x_n=r\sin \phi$) is given by		
		 $\{\phi_1\le \phi \le \phi_2, 0\le r\le 1\}$, for some angles $\phi_1, \phi_2\in (0,\pi)$. Let the projection of $W_i$ onto the plane $x_1x_n$ be the point $\tau=(\cos\phi_0, 0,...,0, \sin\phi_0)$, where $\phi_1\le \phi_0\le \phi_2$, see Figure \ref{Fig:icecream-corners}.
		 
		 For parameters $\lambda $ and $\mu$ such that $\phi_1\le \lambda \le \mu \le  \phi_2$, define $$Q_{\lambda,\mu}=\{x\in\mathbb R^n:   x_1= r\cos\phi, x_n=r\sin\phi, r\ge 0, \lambda\le \phi\le\mu\}.$$  
		Recall  $  U_i=\{ t \xi:   t \ge 0, \ \xi \in \Omega_i \}$.  

\begin{figure}[h!]
 \centering
\begin{tikzpicture}
\node[anchor=south west, inner sep=0] (img) at (0,0)
        {\includegraphics[scale=0.23]{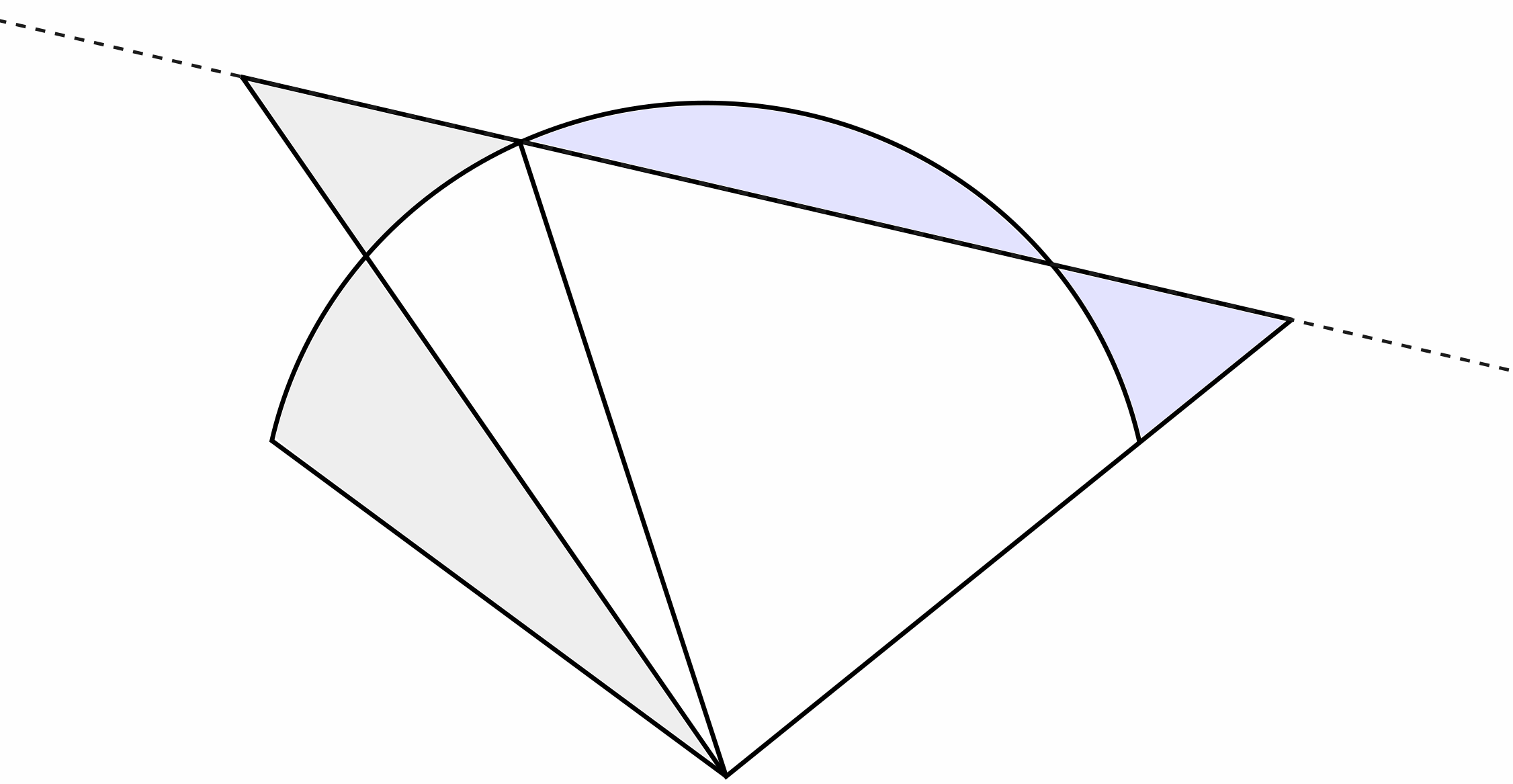}};
\draw [fill=black] (6.97,3.5) circle (1.5pt);
\draw [fill=black] (5.475,3.5) circle (1.5pt);
\node[above right] at (13.8,3.4) {$\psi_0$};
\node[above right] at (12.6,4.1) {$\psi_0^+$};
\node[above right] at (12.6,3.1) {$\psi_0^-$};
\node[above right] at (7.05,3.3) {$c(\hat\Omega_i)$};
\node[above left] at (5.4,3.4) {$c(\Omega_i)$};
\node[above right] at (10.3,3.7) {$I$};
\node[above right] at (6.34,5.5) {$II$};
\node[above right] at (3,5.3) {$III$};
\node[above right] at (3,3) {$IV$};
\draw[thick, dashed] (5.475,3.5) -- (4.2,7.397);
\node[above right] at (3.85,7.4) {$v_i^\perp$};
\draw[thick, dashed] (6.6,0.1) -- (7.39,7.359);
\node[above right] at (7.2,7.35) {$w_i^\perp$};
\end{tikzpicture}
\caption{Construction of cone $\hat \Omega_i$ from $\Omega_i$ in Step 5.}
	\label{Fig:icecream}
\end{figure}
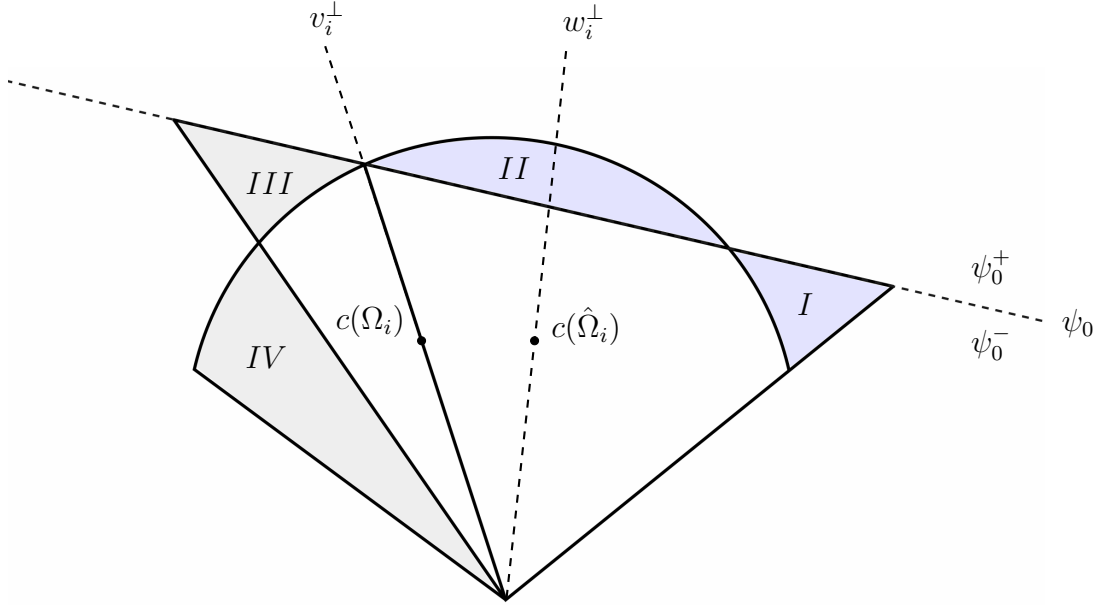

        Now we find a specific hyperplane $\psi_0 \in \Psi_i$ and a parameter $\mu_0\in(\phi_0, \phi_2)$ that are used to define a cone $\hat \Omega_i$. 
		For $\psi\in \Psi_i$, let $\psi^-$ be the half-space bounded by $\psi$ that contains the origin, and $\psi^+$ be its complement. 
		If $\psi$ is tangent to the cylinder $\{x\in \mathbb R^n: x_1^2+x_n^2= 1\}$, then $$\mathrm{vol}_n(\Omega_i\cap Q_{\phi_1,\phi_0}\cap \psi^+)=0 \qquad\mbox{and}\qquad\mathrm{vol}_n((U_i\setminus \Omega_i)\cap Q_{\phi_1,\phi_0}\cap \psi^-)>0.$$ 
		On the other hand, if $\psi$ contains the point $(\cos \phi_1,0,...,0,\sin\phi_1)$, then 
		$$\mathrm{vol}_n(\Omega_i\cap Q_{\phi_1,\phi_0}\cap \psi^+)>0 \qquad\mbox{and}\qquad\mathrm{vol}_n((U_i\setminus \Omega_i)\cap Q_{\phi_1,\phi_0}\cap \psi^-)=0.$$ 
		 By continuity and the intermediate value theorem, there is $\psi_0\in \Psi_i$ such that 
		$$\mathrm{vol}_n(\Omega_i\cap Q_{\phi_1,\phi_0}\cap \psi_0^+) = \mathrm{vol}_n((U_i\setminus \Omega_i)\cap Q_{\phi_1,\phi_0}\cap \psi_0^-).$$
In Figure \ref{Fig:icecream}, the sets  $(U_i\setminus \Omega_i)\cap Q_{\phi_1,\phi_0}\cap \psi_0^-$ and $\Omega_i\cap Q_{\phi_1,\phi_0}\cap \psi_0^+$  are depicted by regions  $I$ and $II$, respectively.

		Now we fix this $\psi_0$ and consider $\mu\in[\phi_0,\phi_2]$. We apply the intermediate value theorem again to find $\mu_0$.	
 If $\mu=\phi_2$, then 
		$$\mathrm{vol}_n(Q_{\phi_0,\mu}\cap (U_i\setminus \Omega_i)\cap \psi_0^-)>0 \qquad\mbox{and}\qquad\mathrm{vol}_n( (\Omega_i\cap Q_{\phi_0,\phi_2}) \setminus Q_{\phi_0,\mu})=0.$$ 
		On the other hand, if $\mu=\phi_0$, then
		$$\mathrm{vol}_n(Q_{\phi_0,\mu}\cap (U_i\setminus \Omega_i)\cap \psi_0^-)=0 \qquad\mbox{and}\qquad\qquad\mathrm{vol}_n( (\Omega_i\cap Q_{\phi_0,\phi_2}) \setminus Q_{\phi_0,\mu})>0.$$
		By the intermediate value theorem, there is $\mu_0\in (\phi_0,\phi_2)$ such that 
		$$\mathrm{vol}_n(Q_{\phi_0,\mu_0}\cap (U_i\setminus \Omega_i)\cap \psi_0^-)=\mathrm{vol}_n( (\Omega_i\cap Q_{\phi_0,\phi_2}) \setminus Q_{\phi_0,\mu_0}).$$
        In Figure \ref{Fig:icecream}, the sets $Q_{\phi_0,\mu_0}\cap (U_i\setminus \Omega_i)\cap \psi_0^-$ and $(\Omega_i\cap Q_{\phi_0,\phi_2}) \setminus Q_{\phi_0,\mu_0}$ are depicted by regions $III$ and $IV$.

		  Thus, we define $$\hat \Omega_i = U_i \cap Q_{\phi_1,\mu_0}\cap  \psi_0^-.$$
		Observe that, by construction,
		$$\mathrm{vol}_n (\hat \Omega_i)=\mathrm{vol}_n (  \Omega_i),$$
		$$\mathrm{vol}_n (\hat \Omega_i\cap v_i^+)=\mathrm{vol}_n (  \Omega_i\cap v_i^+),$$
		and the centroid of $\hat \Omega_i$ lies in $v_i^+$.  Assume that $c(\hat\Omega_i) \in w_i^\perp$, where $w_i^\perp = \{x\in \mathbb R^n: x_1=r\cos \delta, x_n =r\sin \delta, r\ge 0\}$, for some angle  $\delta\in (\phi_1, \phi_0]$, and $\langle p, w_i\rangle$ has the same sign as $\langle p, v_i\rangle$.
		Hence,
		\begin{equation}\label{Omega-to-Omega-hat}\frac{\mathrm{vol}_n (\Omega_i\cap v_i^+)}{\mathrm{vol}_n (\Omega_i)} = \frac{\mathrm{vol}_n (\hat\Omega_i\cap v_i^+)}{\mathrm{vol}_n (\hat\Omega_i)}\ge \frac{\mathrm{vol}_n (\hat\Omega_i\cap w_i^+)}{\mathrm{vol}_n (\hat\Omega_i)}.
        \end{equation}

		Note that $\hat \Omega_i$ is a cone that is the convex hull of the origin and an $(n-1)$-dimensional convex set $A_i$  in $\psi_0$. Since $$\mathrm{vol}_n (\hat\Omega_i)=\frac1n h\, \mathrm{vol}_{n-1} (A_i)\quad\mbox{and}\quad \mathrm{vol}_n (\hat\Omega_i\cap w_i^+)=\frac1n h\, \mathrm{vol}_{n-1} (A_i\cap w_i^+),$$
		where $h$ is the distance from the origin to the hyperplane $\psi_0$, and 
		 since the centroid of $A_i$ lies on the same ray from the origin as the centroid of $\hat \Omega_i$, so it also lies in $w_i^\perp$,
		we obtain 
		\begin{equation}
		\label{Omega-hat-to-A}\frac{\mathrm{vol}_n (\hat\Omega_i\cap w_i^+)}{\mathrm{vol}_n (\hat\Omega_i)}=\frac{\mathrm{vol}_{n-1} (A_i\cap w_i^+)}{\mathrm{vol}_{n-1} (A_i)}\ge \left(1-\frac{1}{n}\right)^{n-1}.
        \end{equation}
		The latter follows from Gr\"unbaum's inequality \eqref{Gru} applied to the $(n-1)$-dimensional set $A_i$. 

        Combining relations \eqref{u-to-v}, \eqref{K-to-Omega}, \eqref{v-to-vi}, \eqref{Omega-to-Omega-hat}, and \eqref{Omega-hat-to-A}, we obtain
		\begin{align*}
\frac{\sigma(K \cap u^+)}{\sigma(K)}&\ge 
\lim_{i\to \infty} \frac{\sigma(K_i \cap v^+)}{\sigma(K_i )}=\lim_{i\to \infty}\frac{\mathrm{vol}_n (\Omega_i\cap v^+)}{\mathrm{vol}_n (\Omega_i)} =\lim_{i\to \infty}\frac{\mathrm{vol}_n (\Omega_i\cap v_i^+)}{\mathrm{vol}_n (\Omega_i)} \\ &=\frac{\mathrm{vol}_n (\hat\Omega_i\cap v_i^+)}{\mathrm{vol}_n (\hat\Omega_i)}\ge \frac{\mathrm{vol}_n (\hat\Omega_i\cap w_i^+)}{\mathrm{vol}_n (\hat\Omega_i)} =\frac{\mathrm{vol}_{n-1} (A_i\cap w_i^+)}{\mathrm{vol}_{n-1} (A_i)}
\ge \left(1-\frac{1}{n}\right)^{n-1},
        \end{align*}
			which concludes the proof in the case when $K_\infty$ is a spherical segment.

		{\bf Step 6.} Now we treat the case where $K_{\infty}$ is a point. We need to show that
		$$
		\lim_{i\to \infty} \frac{\sigma(K_i \cap u^+)}{\sigma(K_i )} \ge \left(1-\frac{1}{n}\right)^{n-1}.
		$$		
		Without loss of generality, we assume that the point $K_\infty$ is located at the north pole $N$. Consider the gnomonic projection centered at $N$. Let $\tilde{K}_i$ be the image of $K_i$ and let $\tilde{q}_i$ be the image of $q_i$, the centroid of $K_i$. 
        
        Since $\{K_i\}$ approaches $N$ in the Hausdorff metric as $i \to \infty$, we can choose a decreasing sequence of positive numbers $\{\varepsilon_i\}$ that converges to zero such that $\tilde{K}_i$ is contained in the ball $\varepsilon_i B_2^{n-1}\subset \{x_n=1\}$ for every $i$.

		Then
		\begin{align*}
			\sigma(K_i)  &= \int_{\tilde   K_i} (1+|x|^2)^{-n/2} \,dx_1\ldots dx_{n-1}\\
			&= \int_{\tilde   K_i} (1 + O(\varepsilon_i^2)) \,dx_1\ldots dx_{n-1}  = (1 + O(\varepsilon_i^2)) \mathrm{vol}_{n-1} (\tilde   K_i).
		\end{align*}
		Similarly,
		\begin{align*}
			\sigma(K_i \cap u^+)  &= \int_{\tilde   K_i \cap u^+} (1+|x|^2)^{-n/2} \,dx_1\ldots dx_{n-1}\\
			&= \int_{\tilde   K_i \cap u^+} (1 + O(\varepsilon_i^2)) \,dx_1\ldots dx_{n-1}  = (1 + O(\varepsilon_i^2)) \mathrm{vol}_{n-1} (\tilde   K_i \cap u^+).
		\end{align*}
		Hence,
		\begin{align}\label{Eqn1}
			\frac{\sigma(K_i \cap u^+)}{\sigma(K_i)} = \frac{\mathrm{vol}_{n-1} (\tilde   K_i \cap u^+)}{ \mathrm{vol}_{n-1} (\tilde   K_i)}(1 + O(\varepsilon_i^2)).
		\end{align}

		Let  $c(\tilde K_i)$ be the centroid of $\tilde K_i$ in the hyperplane $\{x_n=1\}$. Note that $c(\tilde K_i)$ does not necessarily lie in $u^\perp$. However, we show that  it is  $O(\varepsilon_i^2)$ close to $u^\perp$.

        Without loss of generality, we assume that $u=e_1$.  For brevity, we write $\nu(x)=(1+|x|^2)^{-\frac{n+1}{2}}$ and $dx=dx_1\ldots dx_{n-1}$. Since $q_i$, and therefore $\tilde q_i$, lie in $u^\perp$, \eqref{proj-of-centroid} yields that
		$$  
        \frac{\int_{\tilde{K}_i} x_1\nu(x) \, dx}{\int_{\tilde{K}_i} \nu(x) \, dx}=\langle \tilde q_i, e_1 \rangle=0.
        $$
        Thus, we obtain
		 \begin{align*}
		 	&|\langle c(\tilde K_i),e_1\rangle |\\
            &= \left |\frac{\int_{\tilde{K}_i} x_1  \, dx}{\int_{\tilde{K}_i}   \, dx } \right|= \left |\frac{\int_{\tilde{K}_i} x_1  \, dx}{\int_{\tilde{K}_i}   \, dx }-\frac{\int_{\tilde{K}_i} x_1\nu(x) \, dx}{\int_{\tilde{K}_i} \nu(x) \, dx }\right|\\
		 	&= \frac{\left | \int_{\tilde{K}_i} \nu(x) \, dx \int_{\tilde{K}_i} x_1  \, dx - \int_{\tilde{K}_i}   \, dx\int_{\tilde{K}_i} x_1\nu(x) \, dx \right|}{\int_{\tilde{K}_i}   \, dx \int_{\tilde{K}_i} \nu(x) \, dx }\\
            &= \frac{\left | \left(\int_{\tilde{K}_i} \nu(x) \, dx \int_{\tilde{K}_i} x_1  \, dx - \int_{\tilde{K}_i}  \, dx  \int_{\tilde{K}_i} x_1  \, dx \right)+\left(\int_{\tilde{K}_i}  \, dx  \int_{\tilde{K}_i} x_1  \, dx-\int_{\tilde{K}_i}   \, dx\int_{\tilde{K}_i} x_1\nu(x) \, dx\right) \right|}{\int_{\tilde{K}_i}   \, dx \int_{\tilde{K}_i} \nu(x) \, dx }  \\
		 	&\le 	 \frac{\left | \int_{\tilde{K}_i} \nu(x) \, dx - \int_{\tilde{K}_i}  \, dx \right| \int_{\tilde{K}_i} |x_1|  \, dx  + \int_{\tilde{K}_i}  \, dx\left| \int_{\tilde{K}_i} x_1  \, dx-  \int_{\tilde{K}_i} x_1\nu(x) \, dx \right|}{\int_{\tilde{K}_i}   \, dx \int_{\tilde{K}_i} \nu(x) \, dx }.	
		 	\end{align*}
		 	Note that 
		$$\left| \int_{\tilde{K}_i} \nu(x) \, dx - \int_{\tilde{K}_i}  \, dx \right|\le O(\varepsilon_i^2) \mathrm{vol}_{n-1} (\tilde   K_i),$$
        and
		$$\int_{\tilde{K}_i} |x_1|  \, dx\le h \,\mathrm{vol}_{n-1} (\tilde   K_i),$$
		where 
		$$h = \max_{x\in \tilde K_i} |\langle x,e_1\rangle|.$$
		Also,
		$$\left| \int_{\tilde{K}_i} x_1  \, dx-  \int_{\tilde{K}_i} x_1\nu(x) \, dx \right|\le h \, O(\varepsilon_i^2) \mathrm{vol}_{n-1} (\tilde   K_i)$$
		and
		$$\int_{\tilde{K}_i} \nu(x) \, dx = (1+O(\varepsilon_i^2))  \mathrm{vol}_{n-1} (\tilde   K_i).$$
		Thus, 
$$	|\langle c(\tilde K_i),e_1\rangle |\le D h\,  \varepsilon_i^2,$$
for some absolute constant $D>0$.

Let $b= \langle c(\tilde K_i),e_1\rangle$. We will estimate the difference 		
\begin{align*}&\left|	\mathrm{vol}_{n-1} (\tilde   K_i \cap \{x_1\ge 0\})  - \mathrm{vol}_{n-1} (\tilde   K_i \cap \{x_1\ge b\})\right|\\ &\qquad\le \int_{ - Dh\varepsilon_i^2}^{D h\varepsilon_i^2}\mathrm{vol}_{n-2}( \tilde  K_i \cap \{x_1=t\})\,  dt\\
	&	\qquad\le   2D h \varepsilon_i^2  \max_{- D h\varepsilon_i^2\le t\le  Dh\varepsilon_i^2}\mathrm{vol}_{n-2}( \tilde  K_i \cap \{x_1=t\}).
	\end{align*}
 Let $t_0\in [- D h\varepsilon_i^2, D h\varepsilon_i^2]$ be a number for which the above maximum is achieved, and let $x_0\in \tilde K_i$ be a point such that $h= |\langle x_0,e_1\rangle|.$ Consider the cone in $\{x_n=1\}$ with apex at $x_0$ and base $ \tilde  K_i \cap \{x_1=t_0\}$.   Let $h_C$ be the height of this cone. Observe that $h_C\ge h-D h\varepsilon_i^2$. 
 Using the fact that the cone is contained in $\tilde K_i$, we obtain
 \begin{align*}&\left|	\mathrm{vol}_{n-1} (\tilde   K_i \cap \{x_1\ge 0\})- \mathrm{vol}_{n-1} (\tilde   K_i \cap \{x_1\ge b\})\right|	\\
 &\qquad \le \frac{1}{n-1} h_C  \mathrm{vol}_{n-2}( \tilde  K_i \cap \{x_1=t_0\}) \frac{2D(n-1)h}{h_C} \varepsilon_i^2\\
 	&\qquad\le \mathrm{vol}_{n-1}(\tilde K_i)  \frac{2D(n-1)}{(1-D\varepsilon_i^2)} \varepsilon_i^2\\
 	&\qquad=  \mathrm{vol}_{n-1}(\tilde K_i)  O(\varepsilon_i^2).
 	\end{align*}
 
 Combining this with \eqref{Eqn1}, we obtain
 		\begin{align}\label{FinalEst}
 	\frac{\sigma(K_i \cap u^+)}{\sigma(K_i)} &= \frac{\mathrm{vol}_{n-1} (\tilde   K_i \cap u^+)}{ \mathrm{vol}_{n-1} (\tilde   K_i)}(1 + O(\varepsilon_i^2)) \\& =\frac{\mathrm{vol}_{n-1} (\tilde   K_i \cap \{x_1\ge b\})}{ \mathrm{vol}_{n-1} (\tilde   K_i)}(1 + O(\varepsilon_i^2))\ge \left(1-\frac{1}{n}\right)^{n-1} \left(1+ O(\varepsilon_i^2)\right), \nonumber
 \end{align}
 where we used Gr\"unbaum's inequality \eqref{Gru} applied to the convex set $\tilde K_i$ in $\mathbb R^{n-1}$, i.e., $\{x_n=1\}$, since the hyperplane $x_1=b$ contains the centroid of $\tilde K_i$.

We recall that, by the construction in Step 1, $\left\{\frac{\sigma(K_i \cap u^+)}{\sigma(K_i)}\right\}_{i = 0}^\infty$ is a non-increasing sequence and $K_0 = K$. Then passing to the limit in \eqref{FinalEst} as $i\to \infty$, we get 
 $$	\frac{\sigma(K \cap u^+)}{\sigma(K)} \ge \lim_{i\to\infty}	\frac{\sigma(K_i \cap u^+)}{\sigma(K_i)} \ge \left(1-\frac{1}{n}\right)^{n-1},$$
 as claimed.

 {\bf Step 7.}       
Consider a sequence $\{K_i\}_{i=1}^\infty$ of spherical simplices, i.e.,  each $K_i$ is a convex body in $\mathbb S^{n-1}_+$  such that its gnomonic projection $\tilde K_i$ is a simplex in the hyperplane $\{x_n=1\}$. We also assume that the centroids of $K_i$ are located at the north pole $N$, and $\tilde K_i\subset \varepsilon_i B_2^{n-1}$, where  $\lim_{i\to\infty}\varepsilon_i=0.$ Additionally, choose $u=e_1$, and let each $\tilde K_i$ have a face that is parallel to the  plane $\{x_1=0, x_n=1\}$ and that lies in the half-space $\{x_1\le 0\}$.

Suppose that the centroid of $\tilde K_i$ lies in the plane $\{x_1=b\}$. 
As was shown in Step~6,
\begin{align*}
\left|	\mathrm{vol}_{n-1} (\tilde   K_i \cap \{x_1\ge 0\})- \mathrm{vol}_{n-1} (\tilde   K_i \cap \{x_1\ge b\})\right|	\le  \mathrm{vol}_{n-1}(\tilde K_i)  O(\varepsilon_i^2).
 \end{align*}

    Since $\tilde K_i$  is a simplex with a facet parallel to the cutting plane $x_1=b$ (and this facet lies in $\{x_1\le 0\}$), we have
    $$\frac{\mathrm{vol}_{n-1} (\tilde   K_i \cap \{x_1\ge b\})}{\mathrm{vol}_{n-1}(\tilde K_i)}=\left(1-\frac{1}{n}\right)^{n-1}.$$
    Thus
    $$\left|	\frac{\mathrm{vol}_{n-1} (\tilde   K_i \cap \{x_1\ge 0\})}{\mathrm{vol}_{n-1}(\tilde K_i)}-\left(1-\frac{1}{n}\right)^{n-1} \right|	=  O(\varepsilon_i^2).
    $$
    Since $$\frac{\sigma(K_i \cap u^+)}{\sigma(K_i)} = \frac{\mathrm{vol}_{n-1} (\tilde   K_i \cap u^+)}{ \mathrm{vol}_{n-1} (\tilde   K_i)}(1 + O(\varepsilon_i^2)),$$
    we obtain 
    $$\lim_{i\to \infty} \frac{\sigma(K_i \cap u^+)}{\sigma(K_i)}= \left(1-\frac{1}{n}\right)^{n-1},$$
    which yields the optimality of the  constant in the statement of Theorem \ref{MainThm}.
    
    \qed
    
{\bf Acknowledgments.} The project was made possible by SQuaRE at the American Institute of Mathematics. The authors thank AIM for providing a supportive and mathematically rich environment.


\begin{thebibliography}{9}

\bibitem[AMMY]{AMMY}D.~Alonso-Gutiérrez, F.~Mar\'in Sola,  J.~Mart\'in Go\~ni,  J.~Yepes Nicol\'as, {\it
A general functional version of Gr\"unbaum's inequality},
J. Math. Anal. Appl. {\bf 44} (2025), no. 1, Paper 129065, 20~pp.

\bibitem[BHLL]{BHLL} B.~Basit, S.~Hoehner, Z.~L\'angi, J.~Ledford,
{\it Steiner symmetrization on the sphere},  Indiana Univ.
Math. J., to appear,  arXiv:2406.10614.



\bibitem[BO]{BO} 
A.~Basu, T.~Oertel,  
{\it Centerpoints: a link between optimization and convex geometry},
SIAM J. Optim. {\bf 27} (2017), no. 2, 866--889.

\bibitem[BHPS]{BHPS} F.~Besau, T.~Hack, P.~Pivovarov, F.~Schuster, {\it Spherical centroid bodies}, Amer. J. Math. {\bf 145} (2023), 515--542.


\bibitem[B]{B} W.~Blaschke, Vorlesungen \"uber Differentialgeometrie. II, Affine Differentialgeometrie,
Springer, Berlin (1923).

\bibitem[BF]{BF} T.~Bonnesen, W.~Fenchel, {Theory of convex bodies}, BCS Associates, Moscow, ID, 1987.

\bibitem[BOS]{BOS} V.-E. Brunel, S. Ohta, J. Serres, {\it A generalization of Gr\"unbaum’s inequality in RCD(0, N)-spaces},  J. Funct. Anal.  {\bf 290} (2026), 111210.




\bibitem[CS]{CS} A.~Cristi, D.~Salas, {\it 
Reducing the large set threshold for Oertel's conjecture on the mixed-integer volume}, Lecture Notes in Comput. Sci., 15620,
Springer, Cham, 2025, 199--212.

\bibitem[FLLMT]{FLLMT} M.~Fradelizi, D.~Langharst, J.~Liu,  F.~Mar\'in Sola, S.~Tang, {\it Gr\"unbaum's inequality for Gaussian and convex probability measures}, J. Math. Pures Appl. {\bf 213} (2026), 103938.

	\bibitem[FMY]{FMY}M.~Fradelizi, M.~Meyer, V.~Yaskin, {\it On the volume of sections of a convex body by cones}, Proc. Amer. Math. Soc. \textbf{145} (2017), no. 7, 3153--3164.

\bibitem[G]{G} {B.~Gr\"unbaum}, {\it Partitions of mass-distributions and of convex bodies by hyperplanes}, Pacific
J. Math. {\bf 10} (1960), 1257--1261.

\bibitem[LY]{LY} B.~Letwin, V.~Yaskin, {\it A generalization of Gr\"unbaum's inequality},  Israel J. Math., to appear, arXiv:2410.04741.

		\bibitem[MNRY]{MNRY}M.~Meyer, F.~Nazarov, D.~Ryabogin, V.~Yaskin,
	{\it Gr\"unbaum-type inequality for log-concave functions}, Bull. Lond. Math. Soc. \textbf{50} (2018), no. 4, 745--752.
    
\bibitem[MSZ]{MSZ} S.~Myroshnychenko, M.~Stephen, N.~Zhang, {\it Gr\"unbaum's inequality for sections},  J. Funct. Anal. {\bf 275} (2018),  no. 9,  2516--2537.

\bibitem[O]{O} T. Oertel, {\it  Integer convex minimization in low dimensions}, Ph.D. thesis, ETH Zurich (2014).

\bibitem[R]{R} A.~Rusciano, {\it A Riemannian corollary of Helly’s theorem}, J. Convex Anal. {\bf 27} (2020), no.4, 1261--1275.
\bibitem[ShY]{ShY} A.~Shyntar, V.~Yaskin, {\em A generalization of Winternitz's theorem and its discrete version},  Proc. Amer. Math. Soc. {\bf 149} (2021), no. 7, 3089-3104. 
					
\bibitem[SY]{SY}	M.~Stephen, V.~Yaskin, {\em Applications of Gr\"unbaum-type inequalities},  Trans. Amer. Math. Soc. {\bf 372} (2019), 6755--6769.
				
\bibitem[SZ]{SZ} {M.~Stephen, N.~Zhang}, {\em Gr\"unbaum's inequality for projections}, J. Funct. Anal. \textbf{272} (2017), no. 6, 2628--2640.
				
\end{thebibliography}
\end{document}